\documentclass[letter,12pt]{article}
\usepackage{amssymb}
\usepackage{amsmath, psfrag, graphics, graphicx}
\usepackage{cite}
\newtheorem{teo}{Theorem}
\newtheorem{lema}[teo]{Lemma}
\newtheorem{prop}[teo]{Proposition}

\newcommand{\R}{\mathbb{R}}
\newcommand{\x}{x}
\newcommand{\esp}{\hspace{0.2cm}}

\newcommand{\N}{\mathbb{N}}

\newcommand{\Z}{\mathbb{Z}}

\newcommand{\dist}{\mathrm{dist}}
\newcommand{\B}{\mathrm{Ball}}

\newcommand{\clo}{\mathbb{T}^1}

\begin{document}

\title{On the dynamics of non-reducible cylindrical vortices}
\author{\Large{Daniel Coronel, \hspace{0.02cm} 
Andr\'es Navas \hspace{0.02cm} \& \hspace{0.02cm} Mario Ponce}}
\maketitle

\date{}

\vspace{-0.2cm}

\noindent{\bf Abstract:} We study the dynamics of Euclidean isometric extensions 
of minimal homeomorphisms of compact metric spaces. Under a general hypothesis of 
homogeneity for the base space, we show that these systems are never minimal, thus 
extending a classical result of Besicovitch concerning cylindrical cascades. Moreover, 
using Anosov-Katok type methods, we construct a topologically transitive isometric 
extension over an irrational rotation with a 2-dimensional fiber.

\vspace{0.2cm}

\noindent{\bf MSC:} 37B05, 37E30, 37F50, 

\vspace{0.8cm}

\noindent{\Large{\bf Introduction}}

\vspace{0.3cm}

This work is a natural companion of \cite{CNP}, where the dynamics of cocycles 
of isometries of nonpositively curved spaces (including $\mathbb{R}^{\ell}$) over 
a minimal dynamics (semigroup action) is studied. While \cite{CNP} is mostly 
devoted to the {\em reducible} case, that is, when there is a continuous invariant 
section (which turns out to be equivalent to the existence of bounded orbits), here 
we concentrate on the opposite ({\em i.e.} non-reducible) case. For simplicity, we 
restrict our attention to actions of $\mathbb{Z}$. Thus, we consider a minimal 
homeomorphism $T: X \to X$ from a compact metric space $X$ to itself, 
and given two continuous functions $\Psi: X \to O(\mathbb{R}^{\ell})$ and 
$\rho: X \to \mathbb{R}^{\ell}$, we consider the dynamics of the 
fibered transformation
$$(x,v) \longrightarrow \big( T(x), \Psi(x)v + \rho(x) \big).$$ 
Such a fibered map will be referred to as a {\em cylindrical vortex}, 
a name that is inspired from that of the classical {\em cylindrical cascades}, which  
correspond to vortices with $\ell = 1$ and $\Psi(x) = \mathrm{Id}$, for all $x \in X$. One  
of the main difficulties of our study is that cylindrical vortices do not 
commute with translations along the fibers. This property holds for 
cylindrical cascades, and it is actually a fundamental tool for the 
study of their dynamical properties ({\em e.g.} the classical proof 
of Gottschalk-Hedlund's theorem \cite{GH}). 

Our main result is a generalization of an old result of Besicovitch 
\cite{besi} to our general framework. Because of technical reasons, we 
restrict ourselves to the case where $X$ is {\em locally homogeneous}, 
that is, for every point $x \in X$ and every neighborhood $V$ of 
$x$, given $y \in V$ there exists a homeomorphism $h_{V,x,y}$ 
sending $x$ into $y$ that is the identity outside $V$. For example, 
topological manifolds and the Cantor set are locally homogeneous. 

\vspace{0.45cm}

\noindent{\bf Main Theorem.} {\em No cylindrical vortex over a 
locally homogeneous space is minimal.}

\vspace{0.45cm}

It is worth to stress that the statement above refers to two-side minimality. 
(All along this work, the word {\em orbit} for an homeomorphism refers 
always to a two-side orbit.) Indeed, the version of this result 
for positive minimality follows from an elementary and 
classical result of Gottschalk; see \cite{G}. 

The validity of our Main Theorem in (fiber) dimension 1 is quite natural; 
for instance, if $X$ is the unit circle, then it follows from and 
important and difficult theorem of Le Calvez and Yoccoz \cite{yoccoz}. 
However, our proof is much simpler and follows the lines of Besicovitch's, 
thought it needs a key modification (we have to consider the case where 
the linear part of the skew dynamics combines $\mathrm{Id}$ and 
$-\mathrm{Id}$.) 

In higher dimension, the situation is rather different, and the arguments 
are of geometric nature. We follow an strategy initiated by Birkhoff 
\cite{birk}, strengthened by P\'erez-Marco \cite{perez-marco}
for germs of 2-dimensional homeomorphisms fixing the origin, and adapted by 
the third-named author for fibered holomorphic maps \cite{ponce}. Basically, 
the idea consists in attaching to each cylindrical vortex a totally 
invariant compact set ``at infinity'', which allows concluding the 
non-minimality. The existence of such a compact set is established 
by an argument of approximation of the base dynamics by periodic ones.

\vspace{0.2cm}

Although non-reducible cylindrical vortices cannot be minimal, they may 
admit minimal invariant closed subsets. This is for example the case if 
there are discrete orbits. However, in the case of a (1-dimensional) cascade, 
the existence of such an orbit gives raise to a nonzero drift. The rest of 
this work deals with zero-drift cylindrical vortices. After slightly 
extending a result due to Matsumoto and Shishikura to this setting 
({\em c.f.} Proposition \ref{ms-general}), we show how subtle is the 
higher-dimensional case. As a ``concrete'' example, we construct a 2-dimensional, 
topologically transitive cylindrical vortex over an irrational rotation of 
the circle such that the angle rotation along the fibers is constant and 
rationally independent of that on the basis. We close this work by discussing 
the arithmetic properties of the pairs of angles that may appear for such an example. 
As a straighforward application of the KAM theory, we show that these pairs must 
satisfy a Liouville type condition provided the corresponding function $\rho$ is smooth.

%%%%%%%%%%%%%%%%%%%%%%%%%%%%%%%%%%%%%%%%%%%%%%%%%%%%%%%%%%%%%%%%%%%%%%%%%%%%%%%%%%%%%%%%%%%%%%%

\section{Non minimality of cylindrical vortices}

\subsection{The 1-dimensional case}
\label{dimension-1}

\hspace{0.45cm}
In this section we prove that 1-dimensional cylindrical vortices cannot be minimal. To do 
this, we need to distinguish three cases: when the linear part of the skew dynamics is the 
identity everywhere, when it coincides with $-\mathrm{Id}$ everywhere, and when it combines 
$\mathrm{Id}$ and $-\mathrm{Id}$. 
The first case was settled by Besicovitch \cite{besi}. The second case follows by 
slightly modifying Besicovitch's arguments. Finally, the third case can be reduced to 
the second one. 

We start by (recalling and) slightly modifying Besicovitch's proof.\footnote{We 
do this in order to avoid the use of the fact that if $F$ is minimal then, 
{\em a-priori}, it must have a dense orbit.} Consider a cylindrical cascade 
\[F \!: (x,v) \mapsto  \big( T(x), v + \rho(x) \big),\] 
and denote by $\Pi$ the projection of $X\times \R$ on $\R$. Obviously, 
if (the $\Pi$-projection of) the orbit of a point is bounded either from 
above or from below, then $F$ cannot be minimal. Assume next that all the 
orbits are unbounded from above and from below. We will show that, in 
this case, all the orbits are proper as maps from $\mathbb{Z}$ into 
$X \times \mathbb{R}$ (compare \S \ref{almost-integrability}). 
Indeed, if the orbit of a point 
$(x,v) \in X \times \mathbb{R}$ is (unbounded from above 
and from below and) not proper, then by examining all 
possible cases one easily convinces that we may choose three 
sequences of integers $n_j < n'_j < n''_j$ such that either 
\begin{equation}\label{eq:enes}
\begin{split}
\lim_{j \to \infty} \Pi (F^{n_j}(x,v)) = 
\lim_{j \to \infty} \Pi (F^{n''_j}(x,v)) \mbox{ belongs to } [-\infty,+\infty),\\ 
\Pi (F^{n'_j}(x,v)) >  j, \quad \textrm{and} \quad
\Pi (F^{n'_j}(x,v)) =  \max_{n_j\le n \le n''_j} \Pi (F^{n}(x,v)). 
\end{split}
\end{equation}
or
\begin{equation}\label{eq:enes2}
\begin{split}
\lim_{j \to \infty} \Pi (F^{n_j}(x,v)) = 
\lim_{j \to \infty} \Pi (F^{n''_j}(x,v)) \mbox{ belongs to } (-\infty,+\infty],\\ 
\Pi (F^{n'_j}(x,v)) <  -j, \quad \textrm{and} \quad
\Pi (F^{n'_j}(x,v)) =  \min_{n_j\le n \le n''_j} \Pi (F^{n}(x,v)). 
\end{split}
\end{equation}
\noindent Let us first consider the case (\ref{eq:enes}).  
Set, for every  $n_j - n'_j \le n \le n''_j - n'_j$,
\[
 z_{j,n} := F^{n'_j+n} \big( x,v - \Pi (F^{n'_j}(x,v)) \big). 
\]
We have: 
\begin{eqnarray}
\label{eq:compacto} z_{j,0}  =  \big( T^{n_j'} (x_0),0 \big), 
&& \textrm{ for every } j\in \N, \\ 
\label{eq:bounded} z_{j,n} = F^{n}(z_{j,0}), \quad \Pi 
(z_{j,n}) \le  0,&& \textrm{ for every } n_j - n'_j \le n \le n''_j - n'_j, 
\quad \textrm{and }\\
\label{eq:Z} n_j - n'_j \rightarrow -\infty, \quad 
n''_j - n'_j \rightarrow +\infty, && \textrm{ as } j\rightarrow +\infty. 
\end{eqnarray}
Due to \eqref{eq:compacto}, passing to a subsequence, we may assume that $z_{j,0}$ 
converges to a point $z_0 := (x_0,0)$. By \eqref{eq:bounded}, $z_{j,n}$ converges 
to $z_n := F^n (z_0)$, for every  $n\in \Z$. Finally, using \eqref{eq:bounded} 
and \eqref{eq:Z}, one may easily see that the orbit $\{z_n\}_{n \in \Z}$ remains 
bounded from above by zero, which contradicts our assumption. Similarly, in 
case (\ref{eq:enes2}), one may conclude the existence of a point whose orbit 
remains bounded from below by zero, which again contradicts our assumption.

\vspace{.2cm}

Consider now a cylindrical vortex of type 
\[F \!: (x,v) \mapsto  \big( T(x), -v + \rho(x) \big).\]
Set $G := F^2$. This is a fibered transformation over the non-necessarily minimal map 
$T^2$. Clearly, if $G$ has a bounded orbit, then the same holds for $F$, which violates 
our hypothesis of minimality for $F$. If not, then the arguments given so far show 
that either:

\noindent -- there is a point $z_0^* = (x_0,v_0)$ whose $G$-orbit is bounded from above, or 

\noindent -- there is a point $z_0^* = (x_0,v_0)$ whose $G$-orbit is bounded from below, or 

\noindent -- all the $G$-orbits are proper.

In the first case, we let 
$$h := \sup_{n \in \mathbb{Z}} \Pi \big( G^n (x_0,v_0) \big).$$ 
With obvious notation, for every $w \in \mathbb{R}$, 
we have 
\begin{eqnarray*}
orb_{F} (x_0,w) &=& orb_{G} (x_0,w) \esp \bigcup \esp F \big( orb_G (x_0,w) \big)\\ 
&\subset& X \times (-\infty, w - v_0 + h] \esp \bigcup \esp F( X \times (-\infty,w - v_0 + h])\\ 
&\subset& X \times \big( -\infty, w - v_0 + h + \| \rho \| \big] \esp \bigcup \esp 
X \times \big[ -w + v_0 - h - \| \rho \|, + \infty \big).
\end{eqnarray*}
Taking $w_0 := -1 + v_0 - h - \| \rho \|$, we see that the $F$-orbit of $(x_0,w_0)$ avoids 
$X \times (-1,1)$. In particular, $F$ is not minimal. The second case can be 
treated similarly. Finally, since $F$ is a proper map, 
if the $G$-orbit of $(x,v_0)$ is proper, then its $F$-orbit 
$orb_{F} (x_0,v_0) = orb_{G} (x_0,v_0) \cup  F \big( orb_G (x_0,v_0) \big)$ 
is also proper, and $F$ cannot be minimal neither. 
 
\vspace{.4cm}

Finally, consider a general 1-dimensional cylindrical vortex 
\[F: (x,v) \mapsto  \big( T(x), \Psi(x) v + \rho(x) \big),\] 
where linear part $\Psi(x)$ can be either $\mathrm{Id}$ or 
$-\mathrm{Id}$ at each 
point.\footnote{For a nice discussion of the measure-theoretical 
cohomological properties of such a map, see \cite{OO}.} Denote 
by $Y \subset X$ the preimage of $\{-\mathrm{Id}\}$ under $\Psi$. This 
is a clopen set. Assuming that it is nonempty, the map $F_Y$ induced 
from $F$ by the first-return map $T_Y$ to the set $Y$ has the form 
\[ F_Y(x,v) = \big( T_Y (x), -v + \tilde\rho(x) \big),\]
where $\tilde \rho \!: Y \to \R$. Since $T$ is a minimal homeomorphism, 
the same must hold for $T_Y$. Therefore, $F_Y$ is a cylindrical vortex 
of the type considered in the second case, thus it cannot be minimal. 
Now, every point in $Y \times \R$ that does not have a dense orbit 
for $F_Y$ also fails to have a dense orbit under $F$, since 
$Y^c \times \R$ is a closed set. Therefore, $F$ is not minimal, 
and this concludes the proof of the Main Theorem for $\ell = 1$.

%%%%%%%%%%%%%%%%%%%%%%%%%%%%%%%%%%%%%%%%%%%%%%%%%%%%%%%%%%%%%%%%%%%%%%%%%%%%%%%%%%%%%%%%%%%%%5

\subsection{The case of higher dimension}
\label{higher-dimension}

\hspace{0.45cm}
In \cite{birk}, to each local planar homeomorphism fixing the origin, Birkhoff associated 
a forward-invariant compact set touching the boundary of the definition domain. In the 
holomorphic case, this construction was refined in \cite{perez-marco} by P\'erez-Marco, 
who constructed a completely invariant compact set touching the boundary. 
In this section, we construct a Birkhoff/P\'erez-Marco like set 
{\em at infinity} (a B/P-M set, for short) for cyclindrical vortices of dimension 
$\ell \geq 2$. To do this, we follow the strategy developed by the third-named author in 
\cite{ponce}, where he extends the construction of P\'erez-Marco to fibered holomorphic 
maps. Next, we use the B/P-M sets to show the non-minimality of these vortices. It is 
worth mentioning that this last idea is not completely new, as Besicovitch's work 
\cite{besi} already includes a remark relating forward non-minimality of 
cylindrical cascades to the existence of Birkhoff invariant sets.  

Given $\ell \geq 2$, let $I \!: \R^{\ell} \to \R^{\ell}$ be an affine Euclidean isometry, that is,
\[ I v = \Psi v + \rho \]
for certain $\Psi \in O (\R^{\ell})$ and $\rho\in \R^{\ell}$. This isometry extends 
continuously to the one-point (Alexandrov) compactification $\R^{\ell} \cup \{\infty\}$ by 
letting $I \infty = \infty$. We will show that the infinity is not an isolated point as 
an invariant object. More precisely, we will show that given an open, bounded set 
$U$, there exists a closed set containing $\infty$ that is completely invariant 
under $I$, touches the boundary $\partial U$, and is contained in $\R^{\ell}\setminus U$. 
We introduce a terminology for this. We say that $K$ is a B/P-M set for $I$ 
avoiding $U$ if the following conditions hold:
  \begin{enumerate}
  \item $K \subset \R^{\ell}\setminus U$.
  \item $K\cup \{\infty\}$ is compact and connected for the one-point 
        compactification topology of $\R^{\ell}$. 
  \item $I(K) = I^{-1}(K)=K$.
  \item $K\cap \partial U \neq \emptyset$. 
  \end{enumerate}

\begin{prop}\label{p1}
For any open, bounded set $U \subset \R^{\ell}$ and any affine isometry 
$I: \R^{\ell} \to \R^{\ell}$, there exists a B/P-M set $K$ for $I$ avoiding $U$.
\end{prop}

To prove this proposition, we need an elementary lemma.
  
\begin{lema}
If the claim of Proposition \ref{p1} holds for $I^2$ and any 
open, bounded set $U \subset \R^{\ell}$, then it also holds for $I$.
\end{lema}

\noindent{\it Proof.} Set $\mathcal{U} := U \cup I(U)$. This is 
an open, bounded set. By hypothesis, there exists a B/P-M set 
$\mathcal{K}$ for $I^2$ avoiding $\mathcal{U}$. Letting 
$K := \mathcal{K} \cup I^{-1}(\mathcal{K})$, it is not 
hard to check that $K$ is a B/P-M set for $I$ avoiding $U$. 
$\hfill\square$
    
\vspace{0.5cm}

\noindent{\it Proof of Proposition \ref{p1}.} We start with the case $\ell=2$. 
Due to the lemma above, we may assume that $\Psi$ preserves orientation, and 
hence we only need to consider two cases:
\begin{itemize}
    \item $\Psi= \mathrm{Id}_{\R^2}$: Take two parallel lines in the direction of 
$\rho$ touching the boundary of the 
open set $U$, and such that $U$ lies in between the band determined by these two lines. Then 
choose $K$ as being the union of the two semi-planes that form the complement of this band. 
(See Figure 1.)
    \item $\Psi$ equals the counterclockwise rotation $R_{\alpha}$ of angle $\alpha \neq 0$: 
The isometry 
$I$ fixes the point $v_0 := (\mathrm{Id} - R_{\alpha})^{-1}\rho$, and $I$ corresponds to the rotation of 
angle $\alpha$ centered at this point. Let $B$ be the smallest open ball centered at this 
point and containing $U$. We then may choose $K := B^c$. (See Figure 2.)
\end{itemize}

%%%%%%%%%%%%%%%%%%%%%%%%%%%%%%%%%%%%%%%%%%%%%%%%%%%%%%%%%%%%%%%%%%%%%%%%%%%%%%%%%%%%%%%%%%%%

\begin{figure*}
\centerline{\includegraphics[height=6cm]{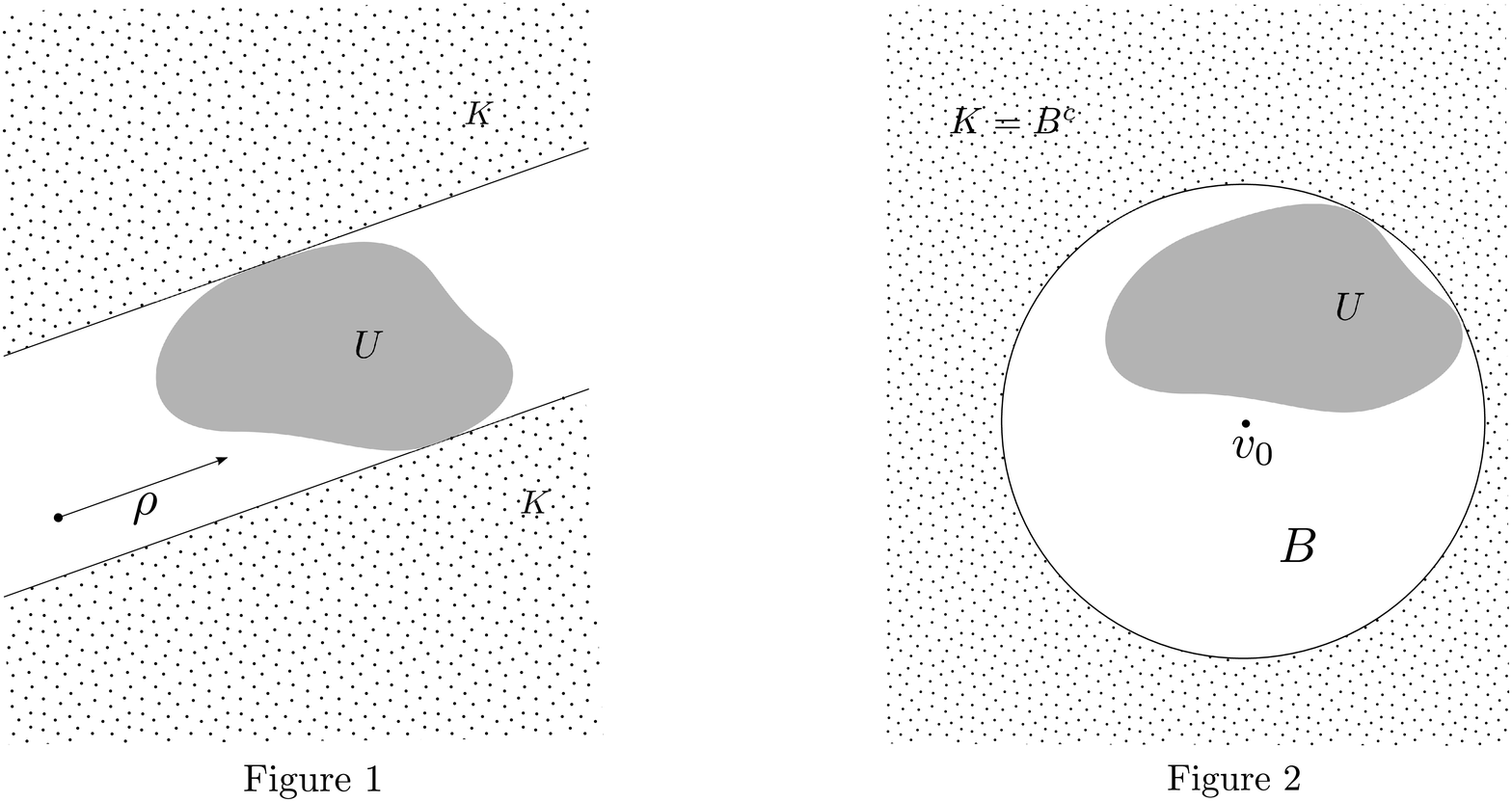}}  
\end{figure*}

%%%%%%%%%%%%%%%%%%%%%%%%%%%%%%%%%%%%%%%%%%%%%%%%%%%%%%%%%%%%%%%%%%%%%%%%%%%%%%%%%%%%%%%%%%%%

When $\ell=3$, in a suitable orthonormal basis, the 
(orientation-preserving) isometry $I$ may written in the form
    \begin{displaymath}
     I v= \left[ \begin{array}{cc}\tilde{\Psi}& \\ 
& 1 \end{array}\right] v+\left[\begin{array}{c}a\\ b 
\\ 0 \end{array}\right]+\left[\begin{array}{c}0\\ 0 \\ c \end{array}\right] 
    \end{displaymath}
for some linear isometry $\tilde \Psi$ of $\R^2$. Let $\tilde U$ be the 
orthogonal projection of $U$ onto $\R^2$, and let $\tilde K$ a B/P-M set 
for $\tilde{I}=\tilde \Psi+[a\ b]^{T}$ avoiding $\tilde U$. (Such a set 
is known to exist because the case $\ell = 2$ has already been settled.) 
Then $K := \tilde K\times \R$ is a B/P-M set avoiding $U$. 

Assume now that the proposition holds for $\ell = k-1$ and $\ell=k$, and consider the case 
$\ell = k+1$. In a suitable orthonormal basis, the isometry may be written in the form
\begin{displaymath}
   I v = \left[ \begin{array}{cc}\tilde{\Psi}& \\ & \hat \Psi \end{array}\right]v 
+ \left[\begin{array}{c}a_1\\ \vdots  \\ a_{k-1}\\ 0\\ 0 \end{array}\right] 
+\left[\begin{array}{c}0\\ \vdots \\ 0\\ a_{k}\\ a_{k+1} \end{array}\right] 
\end{displaymath} 
where $\tilde \Psi, \hat \Psi$ are linear isometries of $\R^{k-1}$ and $\R^{2}$, respectively. 
Let $\tilde U$ be the orthogonal projection of $U$ onto $\R^{k-1}$, and let $\tilde K$ be a   
B/P-M invariant for $\tilde{I}=\tilde \Psi + [a_1 \ldots a_{k-1}]^{T}$ avoiding $\tilde U$. 
Then $K := \tilde K\times \R^2$ is a B/P-M set for $I$ avoiding $U$. $\hfill\square$
    
\vspace{0.35cm}

In what follows, we give an analog of the preceding proposition for fibered isometries. 
More precisely, we consider a cylindrical vortex $F \!: X \times \R^{\ell} \to X \times \R^{\ell}$ 
over a minimal homeomorphism $T \!: X \to X$. We say that an open, bounded set 
$U \subset X \times \R^{\ell}$ is a {\em tube} for $F$ if for every $x\in X$, the 
fiber $U_{x}$ is nonempty. A closed set $K$ is a B/P-M set for $F$ avoiding 
$U$ if the following conditions hold:
\begin{enumerate}
  \item $K\subset \left(X\times \R^{\ell}\right)\setminus U$.
  \item For each $x\in X$, the fiber $K_x \cup \{\infty\}$ is compact and connected 
for the one-point compactification topology of the fiber $\R^{\ell}$. 
  \item $F(K) = F^{-1}(K) = K$.
  \item There exists $x^*\in X$ such that $K_{x^*} \cap \partial U_{x^*} \neq \emptyset$. 
\end{enumerate}

We first consider a fibered dynamics over the finite 
space $X := \Z_p = \Z / p\Z $ with the basis homeomorphism 
$T(j)=j+1$. In other words, given a finite family $\{ I_0, I_1, \dots, I_{n-1} \}$ 
of affine isometries of $\R^{\ell}$, with $\ell \geq 2$, we consider the cylindrical vortex
    \begin{eqnarray*}
    F \!: \Z_p \times \R^{\ell}&\longrightarrow& \Z_p \times \R^{\ell}\\
    (j, v)&\longmapsto&\left(j+1, I_j v \right).
    \end{eqnarray*} 

\begin{prop}\label{ewpp1}
For every tube $U \subset X \times \mathbb{R}^{\ell}$, there 
exists a B/P-M invariant set $K$ for $F$ avoiding $U$.
\end{prop}

\noindent{\it Proof.} For simplicity, we only deal with the case 
$p = 2$, leaving the (slightly more elaborate) general case to the 
reader. Let $U_0,U_1$ be the two (nonempty) open, bounded subsets of 
$\mathbb{R}^{\ell}$ such that $U = (\{ 0 \} \times U_0) \bigcup (\{ 1 \} \times U_1)$, 
and let $\mathcal{U} := U_0 \cup I_0^{-1}(U_1)$. By Proposition \ref{p1}, there 
exists a B/P-M invariant set $K_0$ for $I_1\circ I_0$ avoiding $\mathcal{U}$. 
Set $K_1 := I_0(K_0)$ and $K := (\{0\} \times K_0) \bigcup (\{1\} \times K_1)$. 
Properties 1. and 2. in the definition above are obvious, whereas 
Property 3. follows from 
\begin{eqnarray*}
    (I_1\circ I_0)^{-1}(K_0)&=&K_0,\\
    I_0^{-1}\left(I_1^{-1}(K_0)\right)&=&K_0,\\
    I_1^{-1}(K_0)&=&I_0(K_0) \esp = \esp K_1.
\end{eqnarray*}
Since $K_0\cap \partial\mathcal{U}\neq \emptyset$, either 
$K_0\cap \partial U_0\neq \emptyset$ or $K_0\cap \partial I_0^{-1}(U_1)\neq \emptyset$. 
Since the second condition is equivalent to $K_1\cap \partial U_1\neq \emptyset$, in 
both cases Property 4. above holds. $\hfill\square$
    
\vspace{0.35cm}

We can now proceed to the general case covered by the Main Theorem. Slightly more generally, 
we will say that $T \!: X \to X$ is {\it approximated by homeomorphisms 
having periodic orbits} if there exists a 
sequence of homeomorphisms $T_n \!: X \to X$ that converges to $T$ uniformly on $X$ so that 
each $T_n$ has a periodic point. Our Main Theorem 
follows directly from the next two propositions.

\begin{prop}
Let $X$ be a locally homogeneous compact metric space. If $T\!: X\to X$ 
is a minimal homeomorphism, then $T$ is approximated by homeomorphisms 
having periodic orbits.
\end{prop}

\noindent{\it Proof.} Given $x\in X$, let $V_n$ be a decreasing sequence 
of neighborhoods of $x$ converging to $\{ x \}$. Let $p_n \in \N$ be 
the first-return time of $x$ to $V_n$ under $T$, and let 
$y_n := T^{p_n} x$. Define $T_n := h_{V_n, x, y_n} \circ T$. Clearly, 
$x$ is periodic for $T_n$ with period $p_n$. Finally, since 
$V_n \to \{x\}$, we have the (uniform) convergence 
$T_n\to T$. $\hfill\square$

\begin{prop}\label{dem-for}
Let $U \subset X \times \R^{\ell}$ be a tube. If $T \!: X \to X$ is approximated 
by homeomorphisms having periodic orbits, then there exists a B/P-M set for 
$F$ avoiding $U$.  
\end{prop}

\noindent{\it Proof.} Let 
$x_n^0, x_n^1 := T(x_n^0), \ldots, x_n^{p(n)-1} := T(x_n^{p(n)-2}) $ 
be a periodic orbit of period $p(n)$ of $T_n$, which we identify with 
$\Z_{p(n)}$. Let $F_n \!: \Z_{p(n)} \times \mathbb{R}^{\ell} \to \Z_{p(n)} 
\times \mathbb{R}^{\ell}$ be the cylindrical vortex defined by 
\[
F_n (x_n^{j},v) = \left(x_n^{j+1}, \Psi(x_n^j)v + \rho(x_n^j)\right).
\] 
Proposition \ref{ewpp1} yields a B/P-M invariant set $K_n$ for $F_n$ avoiding 
$U$. Let $\hat K_{n} \subset X \times \overline{\R}^{\ell}$ be the compact 
set resulting from $K_n$ by attaching the section $X \times \{\infty\}$ 
to it, that is, $\hat{K}_n := K_n \cup (X \times \{\infty\})$. Taking 
an appropriate subsequence, we may assume that there exists a connected 
compact set $K \subset X \times \overline{\R}^{\ell}$ that is the limit 
of $K_n$ for the Hausdorff topology on compact sets. Since $\partial U$ is 
a compact set and, for each $n\in \N$, one has $K_n\cap \partial U\neq \emptyset$, 
the intersection $K \cap \partial U$ is nonempty. We denote $\hat K := 
K\cup (X\times\{\infty\})$. Since $T_n \to T$ uniformly, 
we have the uniform convergence $F_n \to F$. Thus, 
$F_n(\hat K_n) \to F (\hat K)$ and $F^{-1}_n (\hat K_n) \to F^{-1}(\hat K)$. 
This implies that $F(K) = F^{-1}(K) = K$, which closes the proof. $\hfill\square$

%%%%%%%%%%%%%%%%%%%%%%%%%%%%%%%%%%%%%%%%%%%%%%%%%%%%%%%%%%%%%%%%%%%%%%%%%%%%%%%%%%%%

\subsection{Two possible generalizations}

\hspace{0.45cm} Recall that the main theorem of \cite{CNP} applies not only to 
cocycles of isometries of a Hilbert space but also to cocycles of isometries of a  
CAT(0) proper space. It is very likely that our Main Theorem here holds in this 
context as well. Indeed, there is a general description of isometries of such 
a space that allows to give an analog of Proposition \ref{p1} for ``higher 
dimensional'' CAT(0)-spaces ({\em e.g.} spaces that are not quasi-isometric 
to the real line). In this situation, this would certainly allow to perform 
a similar procedure to get a B/P-M set and show the non-minimality, whereas for 
the ``1-dimensional case'', it should not be very difficult to adapt the arguments 
of \S \ref{dimension-1}. We do not carry out the details of all of this here 
because we do not see any interesting application except for spaces for which 
the arguments apply without major modifications ({\em e.g.} hyperbolic spaces 
$\mathbb{H}^n$). In this direction, it is worth pointing out that, though 
no cylindrical vortex of isometries of $\mathbb{H}^2$ over an irrational 
rotation is minimal, the action on the product of the circle and the 
boundary of the Poincar\'e disk appears to be minimal in many 
cases \cite{avila-krikorian}.

\vspace{0.1cm}

Perhaps more interesting is trying to settle the infinite dimensional case of the 
Main Theorem. Indeed, most of the arguments of \S \ref{higher-dimension} strongly 
depend on the fact that the fiber space, namely $\R^{\ell}$, is locally compact 
(the same would apply to proper CAT(0)-spaces). This turns natural the following

\vspace{0.4cm}

\noindent{\bf Question.} Does there exist a {\em minimal} cylindrical vortex with 
infinite-dimensional fiber~? 

\vspace{0.34cm}

Recall that, by \cite[Exercise 5.3.15]{book}, no isometry of a Hilbert space 
can be minimal (it cannot be even topologically transitive ``at large scale''). 
However, isometries without fixed points but having recurrent points do exist: 
see \cite{c-t-v} and \cite[Exercise 5.2.26]{book}. The situation might be 
compared with that of general linear maps: in finite dimension, topological 
transitivity is impossible, while in infinite dimension, even topological 
mixing may hold \cite{book-linear}.

%%%%%%%%%%%%%%%%%%%%%%%%%%%%%%%%%%%%%%%%%%%%%%%%%%%%%%%%%%%%%%%%%%%%%%%%%%%%%%%%%%%%%%%%%%%

\section{Minimal invariant sets}

\subsection{Almost-reducibility and proper orbits}
\label{almost-integrability}

\hspace{0.45cm}
The next proposition is folklore but difficult to find in the literature  
stated in this way (compare \cite{krieger,OP}); we include the proof just 
for the convenience of the reader. For the statement, notice that for a 
cylindrical cascade $F \! : (x,v) \to \big( T(x), v + \rho(x) \big)$ 
and each $n \in \mathbb{N}$,
$$F^n (x,v) = \big( T^n (x), v + \rho_n(x) \big),$$
where $\rho_n$ denotes the Birkhoff sum 
$$\rho_n (x) := \sum_{j=0}^{n-1} \rho \big( T^j(x) \big).$$

\begin{prop} \label{caso-cascadas}
The following properties are equivalent:
\begin{enumerate}
\item There exists a family of continuous sections $\varphi_k \!: X \to \mathbb{R}$ 
that is almost invariant under the skew action. In other words, the associate 
cohomological equation can be solved in reduced cohomology:  
$$\lim_{k \to +\infty} \left( \sup_{x \in X} 
\big| \rho(x) - \big[ \varphi_k (T(x)) - \varphi_k (x) \big] \big| \right) = 0.$$
\item We have the convergence  
$$\lim_{n \to +\infty} \left( \sup_{x \in X} 
\Big| \frac{\rho_n (x)}{n} \Big| \right) = 0.$$ 
\item The union of proper orbits has zero $\mu$-measure 
for every $T$-invariant probability measure $\mu$.
\item For every $T$-invariant (ergodic) probability measure $\mu$ on $X$,
$$\int_{X} \rho(x) \hspace{0.015cm} d \mu(x) = 0.$$
\end{enumerate}
\end{prop}

\noindent{\it Proof.} We give the proof of more implications than necessary 
because they will useful for the further discussion of general cylindrical vortices.

\noindent ${\it 1. \to 2.}$ Given $\varepsilon > 0$, let $k = k(\varepsilon)$ 
be such that for all $x \in X$,
$$\big| \rho(x) - [ \varphi_k (T(x)) - \varphi_k (x) \big| \leq \varepsilon.$$
For each $n \in \mathbb{N}$,
$$\big| \rho_n (x) - [ \varphi_k (T^n(x)) - \varphi_k (x) ] \big| 
\leq n \varepsilon,$$
thus
$$\Big| \frac{\rho_n (x)}{n} \Big| 
\leq \frac{2 \| \varphi_k \|_{\infty}}{n} + \varepsilon.$$
Passing to the limit, this yields, with an uniform rate, 
$$\lim_{n \to + \infty} \Big|\frac{\rho_n(x)}{n}\Big| \leq \varepsilon.$$
Since this holds for all $\varepsilon > 0$, we have the desired uniform 
convergence to zero.

\noindent ${\it 2. \to 1.}$ For each $k \in \mathbb{N}$, let 
$$\varphi_k (x) : = - \frac{\rho_1 (x) + \rho_2 (x) + \cdots + \rho_k(x)}{k}.$$
We have 
$$\varphi_k (T(x)) - \varphi_k(x) = \frac{1}{k} \sum_{j=1}^{k} 
\big[ \rho_j(x) - \rho_j(T(x)) \big] = 
\frac{1}{k} \sum_{j=1}^{k} \big[ \rho(x) - \rho (T^{j}(x)) \big] 
= \rho(x) - \frac{\rho_k(T(x))}{k},$$
and the last expression converges uniformly to $\rho(x)$.

\noindent ${\it 1. \to 4.}$ For every $T$-invariant probability measure 
$\mu$ and each $k \in \mathbb{N}$, we have
$$\int_X \big[ \varphi_k (T(x)) - \varphi(x) \big] d \mu = 0.$$
Thus,
$$\int_X \rho \hspace{0.05cm} d \mu 
= \int_X \lim_k \big[ \varphi_k \circ T - \varphi \big] d \mu 
= \lim_k \int_X \big[ \varphi_k (T(x)) - \varphi(x) \big] d \mu = 0.$$

\noindent ${\it 4. \to 1.}$ Let $\mathbb{L}$ be the closure of the subspace of 
$C(X)$ spanned by the functions of the form $\varphi \circ T - \varphi$, where 
$\varphi \in C(X)$. If $\rho$ does not belong to $\mathbb{L}$, the Hahn-Banach 
separation theorem provides us with a linear functional $I$ that restricted to 
$\mathbb{L}$ is zero and $I (\rho) = 1$. Such an $I$ comes from integration 
with respect to a signed probability measure on $X$. Since $I(\mathbb{L}) = \{0\}$, 
this measure is $T$-invariant. Finally, the Hahn decomposition theorem yields 
an invariant probability measure for which the integral of $\rho$ is nonzero.

\noindent ${\it 2. \to 3.}$ This implication directly follows from a classical 
lemma due to Kesten \cite{kesten} (it can be also derived from a well-known 
lemma of Atkinson \cite{atkinson-lema}).

\noindent ${\it 3. \to 4.}$ This follows directly from the Birkhoff ergodic theorem. 
$\hfill\square$

\vspace{0.4cm} 

Notice that $F$ is reducible if and only if the cohomological equation 
$$\rho (x) = \varphi (T(x)) - \varphi (x)$$
has a continuous solution $\varphi \!: X \to \mathbb{R}$. We will hence 
say that the cylindrical cascade $F$ is {\em almost-reducible} if the 
equivalent conditions of the preceding proposition hold. 

The situation for a cylindrical vortex 
$$F \!: (x,v) \to \big( T(x), \Psi(x) v + \rho(x) \big)$$ 
is less transparent. First, in order to introduce a drift-like condition, notice that 
if we define $\rho_n \!: X \to \mathbb{R}$ (and $\Psi_n \!: X \to O (\mathbb{R}^{\ell})$) by 
$$F^n (x,v) =: \big( T^n(x), \Psi_n (x) v + \rho_n(x) \big),$$
then for all $m,n$ in $\mathbb{Z}$, we have
$$\rho_{m+n}(x) = \Psi_n (T^m(x)) \big( \rho_m(x) \big) + \rho_n (T^m (x)).$$
In particular, 
$$\big\| \rho_{m+n}(x) \big\| \leq \big\| \rho_m(x) \big\| 
+ \big\| \rho_n (T^m (x)) \big\|.$$
By the sub-additive ergodic theorem \cite{krengel}, 
for every $T$-invariant ergodic probability measure $\mu$, the value of 
$$\frac{\|\rho_n (x)\|}{n}$$
converges to a limit (drift) $D = D(\mu)$ for $\mu$-almost 
every point $x \in X$. 

The main difference here is that Proposition \ref{caso-cascadas} does not extend 
to cylindrical vortices, even to those with $\Psi \equiv \mathrm{Id}$. More precisely, 
the equivalence between conditions 1. and 2. still holds with an analogous (direct) 
proof. (See \cite{BN} for a more general version of this fact.) 
Nevertheless, all the arguments relying on ergodic type theorems fail. 
As a matter of example, let us consider Yoccoz' example from \cite{yoccoz}. 
This is an irrational rotation of the 2-torus 
$T \!: (x,y) \mapsto (x + \alpha,y + \beta)$ together with two continuous 
functions $\hat{\rho} = \hat{\rho} (x)$ and $\check{\rho} = \check{\rho} (y)$, 
both of zero integral, such that for almost every $(x,y) \in \mathbb{T}^2$,  
$$\lim_{n \to \pm \infty} \big[ |\hat{\rho}_n (x)| + |\check{\rho}_n (y)| \big] = \infty.$$ 
Then letting $\Psi \equiv \mathrm{Id}$ and $\rho := (\hat{\rho},\check{\rho})$, almost 
every orbit of the induced cylindrical vortex $F$ on $\mathbb{T}^2 \times \mathbb{R}^2$ 
is proper. However, for {\em every} point $(x,y) \in \mathbb{T}^2$, we have 
$$\limsup_{n \to + \infty} \frac{\| \rho_n (x,y) \|}{n} = 
\limsup_{n \to + \infty} \frac{|\hat{\rho}_n (x)| + |\check{\rho}_n(y)|}{n} = 0,$$
where the convergence is uniform in $(x,y)$. In particular, $D = D(\mathrm{Leb}) = 0$.

%%%%%%%%%%%%%%%%%%%%%%%%%%%%%%%%%%%%%%%%%%%%%%%%%%%%%%%%%%%%%%%%%%%%%%%%%%%%%%%%%%%%%%%%%%%

\subsection{On a theorem of Matsumoto and Shishikura}

\hspace{0.45cm}
One of the major interests on proper orbits is that they are nontrivial minimal 
invariant closed sets. It easily follows from Denjoy-Koksma's inequality that for any 
almost-reducible cylindrical cascade over an irrational circle rotation, there is no 
such orbit provided that the function $\rho$ has finite total variation (without the 
last assumption, proper orbits may appear; see \cite{sch} and \cite{besi}; see also 
\cite{KS} for recent simpler examples). Actually, as it is shown by Matsumoto and 
Shishikura in \cite{MS}, no nonempty, proper, minimal invariant closed set can 
appear in this situation. A slight extension of this result is our next 

\vspace{0.1cm}

\begin{prop} \label{ms-general}
Let $F$ be an almost-reducible, 1-dimensional cylindrical vortex over 
an irrational rotation of the circle. If the corresponding function $\rho$ has finite 
total variation, then $F$ admits no nonempty, proper, minimal invariant closed set. 
\end{prop}

\noindent{\it Proof.} Let $F \!: (x,v) \to \big( x + \alpha, -v + \rho(x) \big)$ 
be a cylindrical vortex satisfying the hypothesis. (If $\Psi \equiv \mathrm{Id}$, then 
Matsumoto-Shishikura's result applies.) Since $F$ is assumed to be almost-reducible, 
the same must hold for $F^2$. Notice that $F^2$ is a cylindrical cascade to which 
Matsumoto-Shishikura's theorem applies. 

Let $C \neq \emptyset$ be a nonempty minimal closed $F$-invariant subset of $\clo \times \R$. 
Given $(x,v) \in C$, denote by $C_{(x,v)}$ the closure of its (full) orbit under $F^2$. We 
will show below that either $C$ is $F^2$-minimal or $C_{(x_0,v_0)}$ is $F^2$-minimal for 
some $(x_0,v_0)$. Before showing this, notice that the previous remark yields either 
$C = X \times \R$ or $C_{(x_0,v_0)} = X \times \R$, respectively. Since 
$C_{(x_0,v_0)} \subset C$, in the second case we still have 
$C = \clo \times \R$, as we wanted to check.

Assume that no $C_{(x,v)}$ is $F^2$-minimal, and let $C^*_{(x,v)} \subsetneq C_{(x,v)}$ 
be a nonempty closed $F^2$-invariant set. Then the set 
$C^*_{(x,v)} \cup F(C^*_{(x,v)})$ is nonempty, closed, $F$-invariant and 
contained in $C$. By minimality, it coincides with $C$. Now choose $(y,w)$ 
(depending on $(x,v)$) in $C_{(x,v)} \setminus C^*_{(x,v)}$. We must have 
$(y,w) \in F(C^*_{(x,v)}) \subset F(C_{(x,v)})$. Thus, the closed set 
$C_{(x,v)} \cap F(C_{(x,v)})$ is nonempty. Since it is $F$-invariant, 
it must coincide with $C$, and by minimality, this easily implies that $C_{(x,v)} = C$. 
Finally, since this holds for every $(x,v) \in C$, this shows that $C$ 
is $F^2$-minimal, which concludes the proof. $\hfill\square$

%%%%%%%%%%%%%%%%%%%%%%%%%%%%%%%%%%%%%%%%%%%%%%%%%%%%%%%%%%%%%%%%%%%%%%%%%%%%%%%%%%%%%%%%%%%%%%%%%

\subsection{An interesting family of cylindrical vortices}
\label{an-example}

\hspace{0.45cm}
Following \cite[Example 3]{CNP}, given two rationally independent 
angles $\alpha,\beta$ and a continuous function $\rho \!: \clo 
\rightarrow \mathbb{C}$, we consider the cylindrical vortex 
$F \!: (x,z) \mapsto \big( x+\alpha, e^{2\pi i\beta}z + \rho(x) \big)$ 
on $\clo \times \mathbb{C}$.

\vspace{0.15cm}

\begin{lema} The map $F$ has zero drift and is conservative 
(that is, it admits no wandering open domain).
\end{lema}

\noindent{\em Proof.} The function $\rho_n$ defined so that 
\begin{equation}
F^n(x,z) = \big(x + n \alpha, e^{2\pi i n \beta} z + \rho_n (x) \big)
\label{iterados}
\end{equation}
may be rewritten in the form
$$\rho_n (x) 
= \sum_{k=0}^{n-1} e^{2 \pi i (n-k-1)\beta} \rho(x + k\alpha) 
= e^{2\pi i n \beta} \sum_{k=0}^{n-1} e^{- 2\pi i (k+1) \beta} \rho(x + k\alpha).$$
Up to the factor $e^{2\pi in \beta}$, this coincides with the $n^{th}$ Birkhoff sum 
$S_n^T (\chi)(x,0)$ at the point $(x,0) \in \mathbb{T}^2$ of the function 
$\chi (x,y) = e^{2\pi i (y-\beta)} \rho(x)$ with respect to the dynamics of 
the rotation $T \!: \mathbb{T}^2 \to \mathbb{T}^2$ of angle 
$(\alpha,-\beta)$. Indeed, 
\begin{equation}\label{suma-de-Birkhoff}
S_n^T (\chi)(x,y) = \sum_{k=0}^{n-1} e^{2\pi i (y - (k+1) \beta)} \rho(x + k\alpha),
\qquad \rho_n (x) = e^{2\pi in \beta} S_n^T (\chi) (x,0).
\end{equation}
In particular, 
\begin{equation}\label{dividido}
\left| \frac{\rho_n(x)}{n} \right| = \left| \frac{S_n^T (\chi)(x,0)}{n} \right|.
\end{equation}
Since $\alpha,\beta$ are rationally independent, the map $T$ is uniquely ergodic. Since 
$\chi$ is continuous, the right side member of (\ref{dividido}) uniformly converges to 
$$\int_{\mathbb{T}^2} \chi(x,y) \hspace{0.04cm} dx \hspace{0.04cm} dy = 
\int_{\clo}\int_{\clo} e^{2\pi i(y-\beta)} \rho(x) 
\hspace{0.04cm} dx \hspace{0.04cm} dy = 
\int_{\clo} e^{2\pi i(y-\beta)} dy \int_{\clo} \rho(x) dx = 0.$$
In particular, the drift of $F$ is zero. 

To show the conservativity of $F$ notice that, due to Atkinson's lemma \cite{atkinson-lema}, 
if we fix $x \in \clo$ and $\varepsilon > 0$, there exist $\bar{x} \in \clo$, 
$y \in \clo$, and $n \in \mathbb{N}$, such that 
$S_n^T (\chi) (\bar{x},y) \leq \varepsilon$ and 
$$\dist (x,\bar{x}) \leq \varepsilon, \quad \dist (\bar{x}, 
\bar{x} + n \alpha) \leq \varepsilon, \quad 
\dist (y,0) \leq \varepsilon, \quad \dist (y,y-n \beta) \leq \varepsilon.$$ 
This implies that both $n\alpha$ and $n\beta$ are 
$\varepsilon$-close to zero. Together with (\ref{iterados}), (\ref{suma-de-Birkhoff}), and 
$$\dist (x,\bar{x}) \leq \varepsilon, \qquad S_n^T (\chi) (\bar{x},y) \leq \varepsilon,$$
this implies that, for any $z \in \mathbb{C}$, the $n^{\mathrm{th}}$-iterate under $F$ 
of the $(\varepsilon,\varepsilon)$-neighborhood of $(x,z)$ is 
$(\varepsilon,\varepsilon(\|z\|+1))$-close to it. Therefore, 
$F$ has no wandering open domain. $\hfill\square$

\vspace{0.5cm}

As shown by the proof above, the dynamics of $F$ is closely related to the 
cylindrical cascade $G$ on $\mathbb{T}^2 \times \mathbb{C}$ defined by 
$$G \big( (x,y), z \big) \mapsto 
\big( (x+\alpha,y-\beta), z + e^{2\pi i(y-\beta)} \rho(x) \big).$$
Besides (\ref{suma-de-Birkhoff}), both maps are related in that $F$ is a 
factor of $G$. Indeed, letting $\Pi \!: \mathbb{T}^2 \times \mathbb{C} 
\to \clo \times \mathbb{C}$ be the {\em proper map} 
defined by 
$$\Pi \big( (x,y), z \big) := (x, e^{-2\pi iy}z),$$
we have
$$F \circ \Pi = \Pi \circ G.$$
These relations allow showing the next

\vspace{0.1cm}

\begin{lema} \label{implica}
The map $F$ is reducible if and only if $G$ is. 
\end{lema}

\noindent{\em Proof.} In view of either the second relation in (\ref{suma-de-Birkhoff}) 
or the fact that $F$ is a factor of $G$ by a proper map, this follows as a direct 
application of the main result of \cite{CNP}. A direct argument proceeds as follows.

Recall that for $G$ being reducible we mean that there exists a 
continuous function $\varphi \!: \mathbb{T}^2 \to \mathbb{C}$ such that, 
for all $(x,y) \in \mathbb{T}^2$,
$$\varphi(x+\alpha,y-\beta) = \varphi(x,y) + e^{2\pi i(y-\beta)} \rho (x).$$
If this holds, then defining the continuous function 
$\varphi_* \! : \clo \to \mathbb{C}$ by 
$$\varphi_* (x) := \int_{\clo} e^{-2\pi iy} \varphi(x,y) dy,$$
we obtain 
\begin{multline*}
\varphi_* (x) 
= \int_{\clo} \big[ e^{-2\pi iy} \varphi(x+\alpha,y-\beta) - e^{-2\pi i\beta} \rho(x) \big] dy\\
= e^{-2\pi i \beta} \int_{\clo} e^{-2\pi i(y-\beta)} \varphi(x+\alpha,y-\beta) dy 
\esp - \esp e^{-2\pi i\beta} \rho(x) 
= e^{-2\pi i \beta} \big[ \varphi_* (x + \alpha) - \rho(x) \big].
\end{multline*}
Hence
$$\varphi_* (x+\alpha) = e^{2\pi i\beta} \varphi_* (x) + \rho(x),$$
which shows that $F$ is reducible. 

Conversely, assume that $F$ is reducible, that is, there exists a 
continuous function $\varphi_* \! : \clo \to \mathbb{C}$ such that 
$$\varphi_* (x+\alpha) = e^{2\pi i\beta} \varphi_* (x) + \rho(x).$$ 
If we multiply by $e^{2\pi i(y-\beta)}$ both sides of this equality, we obtain 
$$e^{2\pi i(y-\beta)} \varphi_* (x + \alpha) 
= e^{2\pi iy} \varphi_* (x) + e^{2\pi i(y-\beta)} \rho(x).$$
Therefore, if we define $\varphi \!: \mathbb{T}^2 \to \mathbb{C}$ 
by $\varphi(x,y) := e^{2\pi iy} \varphi_* (x)$, we have 
$$\varphi(x+\alpha,y-\beta) = \varphi(x,y) + e^{2\pi i(y-\beta)} \rho (x),$$
thus showing that $G$ is reducible. $\hfill\square$

\vspace{0.6cm}

A more elaborate relation between $F$ and $G$ is given by the next

\vspace{0.1cm}

\begin{prop}\label{trans-top}
The map $F$ is topologically transitive if and only if $G$ is. 
\end{prop}

\vspace{0.1cm}

To show (the difficult implication of) this proposition, we will strongly use a deep 
theorem due to Atkinson that characterizes the failure of topological transitivity 
by the existence of reducible linear factors \cite{atkinson}. In our case, 
this may be stated as follows:

\vspace{0.15cm}

\begin{teo} {\em {\bf [Atkinson]}} Assuming that $G$ is conservative 
and non-reducible, a necessary and sufficient condition for its topological 
transitivity is that there is no $\theta \in \clo$ such that the 
1-dimensional cylindrical cascade
\begin{equation}\label{factor}
\big( (x,y), t \big) \mapsto 
\big( (x+\alpha,y-\beta), t + \big\langle e^{2\pi i\theta}, e^{2\pi i(y-\beta)} \rho(x) \big\rangle \big) 
\end{equation}
is reducible, where $\langle \cdot , \cdot \rangle$ stands 
for the inner product of vectors in $\mathbb{R}^2 \sim \mathbb{C}$.
\end{teo}

\vspace{0.4cm}

\noindent{\em Proof of Proposition \ref{trans-top}.} 
Since $F$ is a factor of $G$, the map $F$ is topologically transitive whenever $G$ is. 

To prove the converse implication, assume first that (\ref{factor}) is reducible 
for some $\theta \in \clo$, that is, there exists a continuous function 
$\varphi \!: \mathbb{T}^2 \rightarrow \mathbb{R}$ 
such that, for all $(x,y) \in \mathbb{T}^2$,
$$\langle e^{2\pi i\theta}, e^{2\pi i(y-\beta)} \rho(x) \rangle = \varphi(x+\alpha,y-\beta) - \varphi(x,y).$$
Then, letting 
$$\varphi_{\vartheta}(x,y) := \varphi(x,y-\vartheta),$$
we have
$$\langle e^{2\pi i\theta'}, e^{2\pi i(y-\beta)} \rho(x) \rangle = 
\varphi_{\theta'-\theta} (x+\alpha,y-\beta) - \varphi_{\theta' - \theta} (x,y).$$
In particular, both cocycles \esp 
$\mathrm{Re} (e^{2\pi i(y-\beta)} \rho(x)) = \langle 1, e^{2\pi i(y-\beta)} \rho(x)\rangle$ \esp 
\esp and \esp $\mathrm{Im} (e^{2\pi i(y-\beta)} \rho(x)) = \langle i, e^{2\pi i(y-\beta)} \rho(x)\rangle$ 
are reducible. Hence, the cocycle $e^{2\pi i(y-\beta)} \rho(x)$ is reducible, which is 
impossible since $G$ is topologically transitive. 

Now, if no cylindrical cascade (\ref{factor}) is reducible, then in order to apply 
Atkinson's theorem for concluding that $G$ is topologically transitive, we need 
to show that $G$ is conservative whenever $F$ is topologically transitive. To do 
this, we first notice that the non-wandering set of $G$ is invariant under both 
$G$ and the translations along the fibers 
$$\big( (x,y), z \big) \mapsto \big( (x,y), z + t \big), \quad t \in \mathbb{C}.$$ 
Therefore, this set is either empty or the whole space $\mathbb{T}^2 \times \mathbb{C}$. 
To exhibit a non-wandering point of $G$ we proceed as follows. Since $F$ 
is topologically transitive, we may chose a point $(x_0,z_0)$ in 
$\clo \times \mathbb{C}$ having a dense orbit under $F$. If 
we denote $(x_n,z_n) := F^n (x_0,z_0)$, this implies in particular that 
there exists a strictly monotone sequence of integers $(n_k)$ such that 
$(x_{n_k},z_{n_k}) \to (x_0,z_0)$ as $k \to \infty$. Since $\Pi$ is a 
proper map, the sequence of subsets 
$\Pi^{-1} (x_{n_k},z_{n_k}) \subset \mathbb{T}^2 \times \mathbb{C}$ 
remains inside a compact set. In particular, there must be a point 
$\big( (\tilde{x}_0,\tilde{y}_0), \tilde{z}_0 \big) \in 
\Pi^{-1} (x_0,z_0)$ such that the sequence 
$\big( (\tilde{x}_{n_k}, \tilde{y}_{n_k}), \tilde{z}_{n_k}) := 
G^{n_k} \big( (\tilde{x}_0,\tilde{y}_0), \tilde{z}_0 \big)$ 
accumulates at some point 
$\big( (\tilde{x}_{\infty},\tilde{y}_{\infty}), \tilde{z}_{\infty} \big)$. 
As it is easy to check, every such an accumulation point is non-wandering 
for $G$. $\hfill\square$

\vspace{0.65cm}

The construction of a map 
$$F \!: (x,z) \mapsto \big( x+\alpha, e^{2\pi i\beta}z + \rho(x) \big)$$ 
that is topologically transitive will be the main issue of the next section. 
Let us close this section by pointing out 
that we do not know whether there exists a cylindrical 
vortex $F$ of the form above that is neither reducible nor 
topologically transitive. This is related to the existence of Yoccoz-like 
examples (see \cite{yoccoz}) associated to functions of a particular 
form, which seems to be a difficult problem. Indeed, the previous 
arguments easily show the following

\vspace{0.15cm}

\begin{prop} If $F$ is neither reducible nor topologically transitive, 
then $G$ is non-conservative. In particular, the function 
$(x,y) \mapsto e^{-2\pi i(y-\beta)}\rho(x)$ does not satisfy 
the Denjoy-Koksma property.
\end{prop}

%%%%%%%%%%%%%%%%%%%%%%%%%%%%%%%%%%%%%%%%%%%%%%%%%%%%%%%%%%%%%%%%%%%%%%%%%%%%%%%%%%%%%%%%%%%%%%%%%

\subsection{A ``concrete'' example}
\label{an-example2} 

\hspace{0.45cm} Quite surprisingly, to perform our construction we will need 
a lemma about the density of a certain set obtained by an arithmetic type 
construction.\footnote{We strongly believe that a much more general result 
should be true; in particular the set $\mathrm{FS}(m,n)$ analogous to that defined 
further one should be dense for all $m > n$. Nevertheless, we were unable 
to produce a conceptual proof of this seemingly interesting fact.} Given 
$q \in \mathbb{N}$, let $t(q) := [\sqrt[3]{q}]$ and $r(q) := [\sqrt{q}]$. 
If $q$ and $t(q)$ are coprime, define $p(q) \in \{1,\ldots,q-1\}$ as the 
inverse (mod $q$) of $t(q)$. Otherwise, set $p(q) := 0.$ Now, let 
$$\mathrm{FS}_q (2,3) := 
\left\{ \Big( \frac{s p(q)}{q}, \frac{s p(q) r(q)}{q} \Big), \esp 
1 \leq s \leq 2 t(q) \right\},$$
where each coordinate is reduced modulo $\mathbb{Z}$ (thus, $\mathrm{FS}_q (2,3)$ 
is though of as a subset of $[0,1]^2$). Finally, let us consider the 
{\em set of frequencies}
$$\mathrm{FS}(2,3):= \bigcup_{q \in \mathbb{N}} \mathrm{FS}_q (2,3).$$

%\vspace{0.1cm}

\begin{lema} \label{F-SET} 
The set $\mathrm{FS}(2,3)$ is dense in $[0,1]^2$.
\end{lema}

\noindent{\em Proof.} For a fixed $m \in \mathbb{N}$, let 
$q := m^6 + 2 m^4 + m^2 + 1 = (m^3 + m)^2 + 1$. One readily 
checks that \esp $t(q) = m^2$ \esp and \esp $r(q) = m^3 + m.$ \esp 
Since \esp $t(q) (m^4 + 2m^2 + 1) = m^6 + 2m^4 + m^2 = q -1,$ \esp 
we have
$$t(q) (q - m^4 - 2m^2 - 1) \equiv 1 \esp (\mathrm{mod} \esp q).$$
Hence,
$$p(q) = q - m^4 - 2m^2 - 1 = m^6 + m^4 - m^2.$$
In particular, modulo $\mathbb{Z}$,
$$\frac{p(q)}{q} = - \frac{m^4 + 2m^2 + 1}{m^6 + 2m^4 + m^2 + 1}.$$ 
Similarly, we have the equality 
$$\frac{p(q) r(q)}{q} = 
- \frac{(m^4 + 2m^2 + 1)(m^3 + m)}{m^6 + 2m^4 + m^2 + 1} 
= -\frac{m^7 + 3m^5 + 3m^3 + m}{m^6 + 2m^4 + m^2 + 1} 
= -m - \frac{m^5 + 2m^3}{m^6 + 2m^4 + m^2 + 1}.$$
Hence, modulo $\mathbb{Z}$,
$$\Big(\frac{s p(q)}{q}, \frac{s p(q) r(q)}{q} \Big) = 
\Big( - \frac{s(m^4 + 2m^2 + 1)}{m^6 + 2m^4 + m^2 + 1}, 
- \frac{s(m^5 + 2m^3)}{m^6 + 2m^4 + m^2 + 1} \Big).$$ 

Let us consider all possible $s$ of the form 
$s (j,k) := jm + k$, where $j,k$ range from $1$ to $m$. Notice that all 
these values satisfy the restriction $1 \leq s \leq 2 t(q) = 2 m^2$, and 
therefore the associated pairs  
$$\Big( - \frac{(jm + k)(m^4 + 2m^2 + 1)}{m^6 + 2m^4 + m^2 + 1}, 
- \frac{(jm + k) (m^5 + 2m^3)}{m^6 + 2m^4 + m^2 + 1} \Big)$$
belong to $\mathrm{FS} (2,3)$. Now, easy computations show that the 
pair above coincides with 
$$\Big( -\frac{j m^5 + km^4 + 2jm^3 + 2km^2 + jm + k}{m^6 + f_4 (m)}, 
-j - \frac{k m^5 + 2km^3 - jm^2 - j}{m^6 + f_4 (m)} \Big)$$
where $f_4$ is a degree-4 polynomial. Since both $j,k$ are 
(positive and) smaller than or equal to $m$, for $m$ large-enough, 
the pair above is very close (modulo $\mathbb{Z}$) to 
$(-\frac{j}{m},-\frac{k}{m}) \sim (1-\frac{j}{m},1-\frac{k}{m})$ 
with an error that converges to zero as $m \to \infty$ (independently 
of $j,k$). As a consequence, every pair of rational numbers in 
$[0,1]^2$ is contained in the closure of $\mathrm{FS} (2,3)$, 
which shows that this set is dense. $\hfill\square$ 

\vspace{0.5cm}

To construct our desired topologically transitive cylindrical vortex 
$$F \!: (x,z) \mapsto \big( x+\alpha, e^{2\pi i\beta}z + \rho(x) \big),$$ 
we will perform a sequence of approximations inspired in 
the classical Anosov-Katok's method \cite{katok}. 
More precisely, we will construct a sequence of skew maps 
$$F_k: (x, z) \mapsto 
\big( x + \alpha_k, e^{2\pi i\beta_k}z + \rho_k (\x) \big), 
\quad \alpha_k \in \mathbb{Q}, \esp \beta_k \in \mathbb{Q},$$
over periodic rotations so that they converge uniformly on compact 
sets. The main point consists in prescribing a sequence of sections 
$\varphi_k \!: \clo \rightarrow \mathbb{C}$ whose images become 
more and more dense in larger and larger regions of $\mathbb{C}$ 
and that are invariant under $F_k$, that is, 
$$\rho_k (\x) 
= \varphi_k (\x + \alpha_k) - e^{2\pi i\beta_k} \varphi_k (\x).$$ 
The construction of the sequence $(\varphi_k)$ is made so that 
\esp $\varphi_k := \varphi_{k-1} + \psi_k,$
\esp where $\psi_k$ has the form 
$$\psi_{k}(\x) := \ell_k \sin(2\pi t_k \x) e^{2\pi i r_{k} \x}.$$
To perform this construction, we will need to define inductively the 
sequences $(\ell_k)$, $(r_k)$, $(t_k)$, $(\alpha_k)$ and $(\beta_k)$. 
The choice is irrelevant for $k = 1$ (we just need to impose the condition 
$\ell_1 \geq 1$). Now, assuming that these values 
have been already defined for $k \in \mathbb{N}$, we let $C_k \geq 1$ be 
a Lipschitz constant for $\varphi_k$ and $D_k \geq 1$ be the supremum 
of the norm of $\varphi_k$. According to Lemma \ref{F-SET}, we 
may choose $(\alpha_{k+1}, \beta_{k+1}) := 
(\frac{s_{k+1} p_{k+1}}{q_{k+1}}, r_{k+1} \alpha_{k+1})$, 
where $p_{k+1}$ is the inverse (modulo $q_{k+1}$) of 
$t_{k+1} := [\sqrt[3]{q_{k+1}}]$ and $r_{k+1} := [\sqrt{q_{k+1}}]$, 
in such a way that the following conditions are satisfied:
\begin{enumerate}
\item \esp $1 \leq s_{k+1} \leq 2 t_{k+1}$.
\item \esp $q_{k+1} > q_k^3$, \esp $q_{k+1}^{5/12} \geq 2^{k+1} q_k$, 
\esp $q_{k+1}^{1/12} \geq 2 k q_{k}^{1/12}$ \esp and \esp 
$q_{k+1} \geq \big[ 2 (k+1) \sum_{j=1}^{k} \ell_j \sqrt{q_j} \big]^3.$
\item \esp $\big| \alpha_{k+1} - \alpha_k \big| \leq 
\frac{1}{2^{k + 2} q_k C_k}$.
\item \esp $\big| \beta_{k+1} - \beta_k \big| \leq 
\frac{1}{2^{k + 2} q_k D_k}$.
\end{enumerate}
Finally, we let $\ell_{k+1} := q_{k+1}^{1/12}$. 

\vspace{0.1cm}

For each $k,n$ in $\mathbb{N}$, 
we define $(\rho_{k})_n \!: \clo \rightarrow \mathbb{C}$ by 
$$(\rho_{k})_n (\x):= \varphi_{k} (\x + n \alpha_k) 
- e^{2\pi in\beta_k} \varphi_k (\x).$$
Then we have 
$$F_k^n (\x, z) = 
\big( \x + n \alpha_k, e^{2\pi in\beta_k} z + (\rho_{k})_n (\x) \big).$$ 
The next lemma yields a quantitative estimate for the convergence 
of the maps $F_k$ as well as some of their iterates.

\vspace{0.1cm} 

\begin{lema} \label{explicit-estimate} 
For each $k \in \mathbb{N}$ and $1 \leq n \leq q_k$, 
\begin{equation}\label{necesaria}
\big| (\rho_{k+1})_n - (\rho_{k})_n \big| \leq 
\frac{\ell_{k+1} q_k s_{k+1}}{q_{k+1}} + 
C_k q_k |\alpha_{k+1} - \alpha_k| + D_k q_k |\beta_{k+1} - \beta_k|.
\end{equation}
\end{lema}

\noindent{\em Proof.} Notice that \esp $(\rho_{k+1})_n (\x) - (\rho_{k})_n (\x)$ \esp equals 
\begin{small}
\begin{eqnarray*} 
&& \hspace{-1cm} 
\varphi_{k+1}(\x + n\alpha_{k+1}) - e^{2\pi in\beta_{k+1}} \varphi_{k+1}(\x) - 
\big[ \varphi_{k}(\x + n \alpha_k) - e^{2\pi i n \beta_k} \varphi_{k}(\x) \big] \\
& = & \varphi_{k+1}(\x + n \alpha_{k+1}) - \varphi_k (\x + n \alpha_k) + 
e^{2\pi i n \beta_k} \varphi_k (\x) - e^{2\pi i n \beta_{k+1}} 
[\varphi_{k}(\x) + \psi_{k+1}(\x)] \\
& = & \big[ \psi_{k+1}(\x + n \alpha_{k+1}) - e^{2\pi i n \beta_{k+1}} \psi_{k+1}(\x) \big] 
+ \big[ \varphi_{k} (\x + n \alpha_{k+1}) - \varphi_k (\x + n\alpha_k) \big] 
+ \varphi_k (\x) \big[ e^{2\pi in\beta_k} - e^{2\pi in\beta_{k+1}} \big].
\end{eqnarray*}
\end{small}Thus, the value of $\big| (\rho_{k+1})_n (\x) - (\rho_{k})_n (\x) \big|$ 
is smaller than or equal to 
\begin{small}
$$\ell_{k+1} \big| \sin ( 2\pi t_{k+1}(\x + n\alpha_{k+1})) 
e^{2\pi i r_{k+1} (\x + n \alpha_{k+1})} - e^{2\pi in\beta_{k+1}} 
\sin(2\pi t_{k+1} \x) e^{2\pi i r_{k+1} \x} \big| + 
C_k n |\alpha_{k+1} - \alpha_k| + D_k n |\beta_{k+1} - \beta_k|.$$
\end{small}The first term of this expression is bounded from above by 
\begin{footnotesize}
\begin{eqnarray*}
&& \hspace{-1.6cm} \big| \ell_{k+1} e^{2\pi ir_{k+1}(\x + n \alpha_{k+1})} 
\big[ \sin(2\pi t_{k+1}(\x + n \alpha_{k+1})) - 
\sin(2\pi t_{k+1} \x) \big] + 
\ell_{k+1} \sin(2\pi t_{k+1} \x) \big[ e^{2\pi i r_{k+1} (\x + n \alpha_{k+1})} 
- e^{2\pi i(n \beta_{k+1} + r_{k+1} \x)} \big] \big|\\ 
\hspace{0.6cm} &\leq&  \ell_{k+1} \big( \{ t_{k+1} n \alpha_{k+1} \} + 
|e^{2\pi inr_{k+1}\alpha_{k+1}} - e^{2\pi in\beta_{k+1}}| \big).
\end{eqnarray*}
\end{footnotesize}Since $p_{k+1} t_{k+1} \equiv 1 \esp (\mathrm{mod} \esp q_{k+1})$, 
we have 
$$\{ t_{k+1} n \alpha_{k+1} \} 
= \left\{ \frac{n s_{k+1} t_{k+1} p_{k+1}}{q_{k+1}} \right\} = 
\left\{ \frac{n s_{k+1}}{q_{k+1}} \right\}.$$
Moreover, since $\beta_{k+1} = r_{k+1} \alpha_{k+1}$, the term 
$|e^{2\pi inr_{k+1}\alpha_{k+1}} - e^{2\pi in\beta_{k+1}}|$ 
vanishes. We thus obtain
$$\big| (\rho_{k+1})_n (\x) - (\rho_{k})_n (\x) \big| \leq 
\frac{\ell_{k+1} n s_{k+1}}{q_{k+1}} + 
C_k n |\alpha_{k+1} - \alpha_k| + D_k n |\beta_{k+1} - \beta_k|,$$
which shows the lemma. $\hfill\square$

\vspace{0.65cm} 

The next lemma deals with the density of the invariant curve 
$\x \mapsto (\x,\varphi_k (\x))$ inside a large  
region of $\clo \times \mathbb{C}$.

\vspace{0.1cm}

\begin{lema} \label{last}
For every $k \ge 2$, the graph of $\varphi_k$ is $(\frac{1}{t_k}, 
\frac{\ell_k t_k}{r_k} + \frac{4\pi}{t_k} \sum_{j = 1}^{k-1} \ell_j r_j)$-dense 
in the cylinder $\clo \times \B \big( 0, \ell_k - 
\sum_{j=1}^{k-1} \ell_j \big)$.
\end{lema}

\noindent{\em Proof.} We first claim that the graph of $\psi_k$ is 
$(\frac{1}{t_k},\frac{\ell_k t_k}{r_k})$-dense 
in $\clo \times \B (0, \ell_k)$. Indeed, given $(\x,z)$ in 
this cylinder, there must exist $\tilde \x \in 
\big[ \x - \frac{1}{2 t_k}, \x + \frac{1}{2 t_k} \big]$ such that 
$\psi_k (\tilde \x)= |z|$. The claim then follows by noticing that 
$\frac{1}{r_k} \leq \frac{1}{t_k}$ and that for every $|s|< \frac{1}{r_k}$,
$$\big| |\psi_k (\tilde \x + s)| - |\psi_k (\tilde \x)| \big| 
\le \ell_k \big( |\sin(2\pi t_k (\tilde \x + s))| - 
    |\sin(2\pi t_k \tilde \x)| \big)
\le \ell_k \frac{t_k}{r_k}.$$

Now, to deal with the graph of $\varphi_k$, we begin by noticing that 
$$|\psi_j'| \leq 2\pi(\ell_j t_j + \ell_j r_j) \leq 4\pi \ell_j r_j.$$
Thus, on each interval $[m/t_k, (m+1)/t_k] \subset \mathbb T^1$, 
the oscillation of $\varphi_{k-1}$ is at most 
$\frac{4\pi}{t_k} \sum_{j = 1}^{k-1} \ell_j r_j.$ 
Since 
$$|\varphi_{k-1}| \leq \sum_{j=1}^{k-1} \ell_j,$$ 
this proves the lemma.
$\hfill\square$

\vspace{0.6cm}

The next lemma deals with the density of a 
certain $F_k$-orbit along the graph of $\varphi_k$.

\vspace{0,1cm}

\begin{lema} \label{densidad-orbita}
The set $\{ F_k^n (0,0) \!: 1 \leq n \leq q_k \}$ is 
$(\frac{1}{t_k}, \frac{\ell_k t_k}{r_k} + 
\frac{4\pi}{t_k} \sum_{j = 1}^{k-1} \ell_j r_j + 
\frac{16\pi}{q_k^{2/3}} \sum_{j=1}^{k} \ell_j r_j)$-dense
in the cylinder $\clo \times \B \big( 0, \ell_k - \sum_{j=1}^{k-1} \ell_j \big)$.
\end{lema}

\noindent{\em Proof.} Since $\varphi_k$ is an invariant section for 
$F_k$, the set of points we are dealing with coincides with 
\begin{equation}\label{set}
\left\{ \left( n \frac{s_k p_k}{q_k}, \varphi_{k} 
\Big( n \frac{s_k p_k}{q_k} \Big) \right) 
\!: 1 \leq n \leq q_k \right\}.
\end{equation}
Since $p_k$ and $q_k$ are coprime and $s_k \leq 2 t_k \leq 2 q_k^{1/3}$, 
the projection on the first coordinate of this set consists of at 
least \esp $q_k^{2/3} / 2$ \esp points uniformly distributed on $\clo$. 
Therefore, the distance between two consecutive points of the set 
(\ref{set}) is less than or equal to 
$$\int_{\frac{2j}{q_k^{2/3}}}^{\frac{2(j+1)}{q_k^{2/3}}} \sqrt{1 + \varphi_k'(\x)^2} 
\esp d \x \esp \leq \esp \frac{4 \max | \varphi_k'|}{q_k^{2/3}} 
\esp \leq \esp \frac{16\pi \sum_{j=1}^{k} \ell_j r_j}{q_k^{2/3}},$$
and the claim of this lemma follows from that of the preceding one. 
$\hfill\square$

\vspace{0.6cm}

To close the construction, notice that for each $1 \leq n \leq q_k$, the properties 
of the inductive construction together with the estimates (\ref{necesaria}) and 
$s_{k+1} \leq 2 q_{k+1}^{1/3} < q_{k+1}^{1/2}$ yield 
\begin{equation}\label{rapidez}
\big| (\rho_{k+1})_n - (\rho_k)_n \big| \leq 
\frac{q_{k+1}^{1/12} q_k q_{k+1}^{1/2}}{q_{k+1}} + 
C_k q_k |\alpha_{k+1} - \alpha_k| + D_k q_k |\beta_{k+1} - \beta_k| 
\leq \frac{q_k}{q_{k+1}^{5/12}} + \frac{1}{2^{k+1}} \leq \frac{1}{2^k}.
\end{equation}
Letting $n := 1$, this shows that $(\rho_k)$ is a Cauchy sequence, 
hence it converges to a continuous function $\rho \!: \clo \rightarrow \mathbb{C}$. 
Moreover, from Property 3. it follows immediately that $(\alpha_k)$ converges to 
some angle $\alpha \in [0,1]$. Similarly, by Property 4, $(\beta_k)$ 
converges to a certain angle $\beta \in [0,1]$. 

Checking that the limit map $F \!: (x,z) \mapsto \big( x+\alpha, e^{2\pi i\beta}z + 
\rho(x) \big)$ is topologically transitive is not very difficult. Indeed, from 
Property 2. it follows that 
$$\frac{1}{t_k} \sum_{j=1}^{k-1} \ell_j r_j \leq \frac{2}{q_k^{1/3}} 
\sum_{j=1}^{k-1} \ell_j q_j^{1/2} \leq \frac{1}{k}.$$ 
Moreover, 
$$\frac{1}{q_k^{2/3}} \sum_{j=1}^{k} \ell_j r_j \leq 
\frac{1}{t_k} \sum_{j=1}^{k-1} \ell_j r_j + \frac{\ell_k r_k}{q_k^{2/3}}
\leq \frac{1}{k} + \frac{\ell_k}{q_k^{1/6}} = \frac{1}{k} + \frac{1}{q_k^{1/12}}.$$
The last two inequalities combined with Lemma \ref{densidad-orbita} imply that 
$\{ F_k^n (0,0) \!: 1 \leq n \leq q_k \}$ is $(\varepsilon_k, \delta_k)$-dense
in the cylinder $\clo \times \B \big( 0, \ell_k - \sum_{j=1}^{k-1} \ell_j \big)$ 
for certain sequences $(\varepsilon_k)$, $(\delta_k)$ converging to zero as $k$ 
goes to infinite. Since $F_k^n (0,0) = \big( n \alpha_k, (\rho_k)_n (0) \big)$, 
using (\ref{rapidez}) and the estimate (valid for $1 \leq n \leq q_k$)
$$\big| n \alpha_k - n \alpha_{k+1} \big| \leq \frac{q_k}{2^{k+2} q_k C_k} 
\leq \frac{1}{2^{k+1}},$$
we conclude that the orbit of $(0,0)$ under the limit map $F$ is dense 
in $\clo \times \B \big( 0, \ell_k - \sum_{j=1}^{k-1} \ell_j \big)$ 
for every $k \in \mathbb{N}$. Finally, since Property 2. yields
$$\sum_{j=1}^{k-1} \ell_j \leq (k-1) \ell_{k-1} = 
(k-1) q_{k-1}^{1/12} \leq \frac{q_k^{1/12}}{2} = \frac{\ell_k}{2},$$ 
we have that \esp 
$\ell_k - \sum_{j=1}^{k-1} \ell_j \geq \frac{\ell_k}{2}$ \esp
goes to infinite together with $k$, thus showing that the $F$-orbit 
of $(0,0)$ is dense in the whole space $\clo \times \mathbb{C}$.

It remains to show that $\alpha$ and $\beta$ are rationally independent. 
Actually, we do not know whether this is always true, but we can ensure 
it provided that the sequence $(q_k)$ satisfies a supplementary condition.

\begin{lema} If the sequence $(q_k)$ satisfies
\begin{itemize}
\item[{\rm 5}. ] $\sum_{j=k}^{\infty}|\alpha_{j+1} - \alpha_j| + 
\sum_{j=k}^{\infty}|\beta_{j+1} - \beta_j| < \frac{1}{k q_k},$
\end{itemize}
then $\alpha$ and $\beta$ are rationally independent.
\end{lema}

\noindent{\em Proof.} Since $p_k \leq q_k$, if we consider the 
representatives of $\alpha_k$ and $\beta_k$ in $[0,1]$, then we have 
$$\alpha_k = \frac{p_k}{q_k}, \quad \beta_k = r_k \alpha_k = 
\frac{r_k p_k}{q_k} - n_k, \quad \mbox{where } n_k \in \mathbb{Z}.$$
Assume that $(p, q, r) \neq (0,0,0)$ is 3-tuple of integers such that 
\esp $p \alpha + q \beta + r = 0.$ \esp On the one hand,
\begin{eqnarray*}
|p \alpha_k + q \beta_k + r|
&=&    |p (\alpha_k - \alpha) + q (\beta_k - \beta) + p \alpha + q \beta + r|\\
&=&    |p (\alpha_k - \alpha) + q (\beta_k - \beta)|\\
&\leq& |p| \sum_{j \geq k} |\alpha_{j+1} - \alpha_j| + 
|q| \sum_{j \geq k} |\beta_{j+1} - \beta_j|.
\end{eqnarray*}
On the other hand, if 
$$p \alpha_k + q \beta_k + r 
= p \!\esp \frac{p_k}{q_k} 
+ q \Big( \frac{r_k p_k}{q_k} - n_k \Big) + r$$ 
equals zero, then  
$$p_k (p + q r_k) = q_k (q n_k - r).$$ 
Since $p_k$ and $q_k$ are coprime, this implies that 
$q_k$ must divide $p + q r_k$. Nevertheless, since $r_k \leq \sqrt{q_k}$
and $(p,q) \neq (0,0)$, this is impossible for a large-enough $k$. 
Therefore, for a large-enough $k \in \mathbb{N}$, the value of 
$p \alpha_k + q \beta_k + r$ is nonzero, and hence
$$|p \alpha_k + q \beta_k + r| 
= \left| \frac{p p_k + q (r_k p_k - n_k q_k) + r q_k}{q_k} \right| 
\geq \frac{1}{q_k}.$$
Therefore, 
$$\frac{1}{q_k} \leq |p| \sum_{j \geq k} |\alpha_{j+1} - \alpha_j| 
+ |q| \sum_{j \geq k} |\beta_{j+1} - \beta_j|,$$
which contradicts Property 5. for $k$ larger than $|p|$ and $|q|$. 
$\hfill\square$

%%%%%%%%%%%%%%%%%%%%%%%%%%%%%%%%%%%%%%%%%%%%%%%%%%%%%%%%%%%%%%%%%%%%%%%%%%%%%%

\section{Reducibility v/s arithmetic properties of the rotation angles}

\hspace{0.45cm} Although the preceding construction provides us with a 
pair of angles $(\alpha, \beta)$ that are rationally independent, it also  
suggests that $\beta$ is very fast approximated by multiples of $\alpha$.  
%(Recall that $\beta_k = r_k \alpha_k$ for all $k \in \mathbb{N}$). 
As we next show, this is the case of every non-reducible smooth cylindrical vortex. 
What follows is inspired (and may be deduced almost directly) from \cite{herman}. 

We say that a pair $(\alpha, \beta) \in \mathbb{T}^1 \times \mathbb{T}^1$ satisfies 
a type-$1$ Diophantine condition, and write $(\alpha, \beta) \in \mathcal{CD}_1$, 
if there exist $C > 0$ and $\tau \geq 0$ such that for every $n \in \mathbb{Z}$,
$$\big| e^{2\pi i(n\alpha-\beta)} - 1 \big| \geq \frac{C}{n^{1+\tau}}.$$
For a fixed irrational $\alpha \in \clo$, denote by $\mathcal{CD}_1^{\alpha}$ the set of 
$\beta \in \clo$ such that $(\alpha, \beta)$ belongs to $\mathcal{CD}_1$. 
Standard arguments show the next

\vspace{0.1cm}

\begin{lema} The set $\mathcal{CD}_1^{\alpha}$ is a countable union 
of closed sets with empty interior and has full Lebesgue measure.
\end{lema}

\vspace{0.1cm}

The next proposition is nothing but a straightforward 
application of the classical baby-KAM theorem.

\begin{prop}
If $\rho:\mathbb{T}\to \mathbb{C}$ is a $C^{\infty}$-function, then 
for every $(\alpha, \beta) \in \mathcal{CD}_1$ the cylindrical vortex
$$F \!:(\x, z)\longmapsto \left(\x+\alpha, e^{2\pi i\beta}z+\rho(\x)\right)$$
is reducible.
\end{prop}

\noindent{\em Proof.} Recall  that reducibility is equivalent to the 
existence of a continuous solution to the cohomological equation
\begin{equation}\label{ceq}
\varphi(\x + \alpha) - e^{2\pi i\beta} \varphi(\x)=\rho(\theta).
\end{equation}
At the level of Fourier series expansion, this is equivalent to that, 
for all $n \in \mathbb{Z}$,
$$\hat \varphi_n = \frac{\hat\rho_n}{e^{2\pi in\alpha} - e^{2\pi i\beta}},$$
where $\hat\varphi_n$ and $\hat\rho_n$ stand for the Fourier coefficients 
of $\phi$ and $\rho$, respectively. Since $\rho$ is a $C^{\infty}$-function, 
$\hat\rho_n$ decreases faster than $1/n^k$ for any $k \geq 1$. The type-1 
Diophantine condition on $(\alpha, \beta)$ allows us to conclude that 
the coefficients $\hat\varphi_n$ defined by the previous equality 
also decrease faster than $1/n^k$ for any $k \geq 1$. Therefore, 
they correspond to the coefficients of a $C^{\infty}$-function 
which solves our cohomological equation. $\hfill\square$

\vspace{0.55cm}

A better result relating the differentiability classes of 
$\rho$ and the solution $\varphi$ can be stated as follows:

\vspace{0.1cm}

\begin{prop} {\bf [Herman]} Let $\rho$ be a  $C^{r}$-function. If 
$(\alpha,\beta) \in \mathcal{CD}_1$ is such that the associated 
$\tau \geq 0$ satisfies $r > \tau + 1$ and $s := r-1-\tau$ is 
not an integer, then the solution $\varphi$ to the equation 
(\ref{ceq}) is a $C^{s}$-function.  
\end{prop}

\vspace{0.1cm}

Finally, concerning the case of ``Liouville pairs'', 
we have the following

\vspace{0.1cm}

\begin{prop}
If $(\alpha, \beta)$ satisfies no type-1 Diophantine condition, then there exists a 
$C^{\infty}$-function $\rho \!: \mathbb{T}^1 \to \mathbb{C}$ such that the equation 
(\ref{ceq}) has no measurable solution. Moreover, $\rho$ can be chosen so that the 
coefficients $\frac{\hat\rho_n}{(e^{2\pi in\alpha} - e^{2\pi i\beta})}$ do not correspond 
to those of any distribution.
\end{prop}

%%%%%%%%%%%%%%%%%%%%%%%%%%%%%%%%%%%%%%%%%%%%%%%%%%%%%%%%%%%%%%%%%%%%%%%%%%%%%%%%%%%%%%%%%%%%%%%%%

\vspace{0.35cm}

\begin{small}

\noindent{\bf Acknowledgments.} We would like to thank to Jan Kiwi for helping 
us in dealing with the sets $\mathrm{FS}(m,n)$, to Raph\"ael Krikorian for 
fruitful discussions on the subject of this work, and to Roman Hric for 
useful information concerning minimal homeomorphisms and for pointing 
out to us the existence of \cite{G}.

\noindent -- D. Coronel was funded by the Fondecyt Post-doctoral Grant 3100092. 

\noindent -- A. Navas was funded by the Math-AMSUD Project DySET.

\noindent -- M. Ponce was funded by the Fondecyt Grant 11090003 and the Math-AMSUD Project DySET.

\end{small}

%%%%%%%%%%%%%%%%%%%%%%%%%%%%%%%%%%%%%%%%%%%%%%%%%%%%%%%%%%%%%%%%%%%%%%%%%%%%%%%%%%%%%%%%%%%%%%%%%

\begin{footnotesize}

%%%%%%%%%%%%%%%%%%%%%%%%%%%%%%%%%%%%%%%%%%%%%%%%%%%%%%%%%%%%%%%%%%%%%%%%%%%%%%%%%%%%%

\vspace{0.25cm}

\noindent{Daniel Coronel}

\noindent{Facultad de Matem\'aticas, PUC}

\noindent{Casilla 306, Santiago 22, Chile}

\noindent{E-mail: acoronel@mat.puc.cl}\\ 

\noindent{Andr\'es Navas}

\noindent{Dpto de Matem\'atica y C.C., USACH}

\noindent{Alameda 3363, Estaci\'on Central, Santiago, Chile}

\noindent{E-mail: andres.navas@usach.cl}\\

\noindent{Mario Ponce}

\noindent{Facultad de Matem\'aticas, PUC}

\noindent{Casilla 306, Santiago 22, Chile}

\noindent{E-mail: mponcea@mat.puc.cl}

\end{footnotesize}


\begin{thebibliography}{dillo83}

\bibitem{katok} {\sc D.V. Anosov \& A.B. Katok.} New examples in smooth ergodic 
theory. Ergodic diffeomorphisms. {\em Trans. Moscow Math. Society} 
{\bf 23} (1970), 1-55.

\bibitem{atkinson-lema} {\sc G. Atkinson.} Recurrence of co-cycles and random walks. 
{\em J. London Math. Soc.} {\bf 13} (1976), 486-488.

\bibitem{atkinson} {\sc G. Atkinson.}  A class of transitive cylinder transformations. 
{\em J. London Math. Soc.} {\bf 17} (1978), 263-270.

\bibitem{avila-krikorian} {\sc A. Avila \& R. Krikorian.} Quasiperiodic 
$\mathrm{SL}(2,\mathbb{R})$ cocycles. In preparation. 

\bibitem{book-linear} {\sc F. Bayart \& \'E. Matheron.} {\em Dynamics of 
Linear Operators.} Cambridge Tracts in Mathematics {\bf 179} (2009).

\bibitem{besi} {\sc A.S. Besicovitch.} A problem on topological transformations of the 
plane II. {\em Proc. Cambridge Philos. Soc.} {\bf 47} (1951), 38-45.

\bibitem{birk} {\sc G. Birkhoff.} Surface transformations and their dynamical 
applications. {\em Acta Math.} {\bf 43}, number {\bf 1} (1922), 1-119.

\bibitem{BN} {\sc J. Bochi \& A. Navas.} A geometric path from zero 
Lyapunov exponents to invariant sections for cocycles. Preprint (2011).

\bibitem{c-t-v} {\sc Y. de Cornulier, R. Tessera, \& A. Valette.} 
Isometric group actions on Banach spaces and representations vanishing 
at infinity. {\em Transform. Groups} {\bf 13} (2008), 125-147.

\bibitem{CNP} {\sc D. Coronel, A. Navas \& M. Ponce.} On bounded cocycles 
of isometries over a minimal dynamics. Preprint (2011), arXiv:1101.3523.

\bibitem{G} {\sc W.H. Gottschalk.} Orbit-closure decompositions and almost periodic 
properties. {\em Bull. of the AMS} {\bf 50} (1944), 915-919

\bibitem{GH} {\sc W.H. Gottschalk \& G.A. Hedlund.} {\em Topological Dynamics.} 
Amer. Math. Soc., Providence, R. I. (1955).

\bibitem{herman} {\sc M. Herman.} Sur les courbes invariantes par les 
diff\'eomorphismes de l'anneau. Vol. 2. {\em Ast\'erisque} {\bf 144} (1986).

\bibitem{kesten} {\sc H. Kesten.} Sums of stationary sequences cannot grow slower 
than linearly. {\em Proc. of the AMS} {\bf 49}, number {\bf 1} (1975), 205-211.

\bibitem{krengel} {\sc U. Krengel.} {\em Ergodic Theorems.} De Gruyter Studies in 
Mathematics {\bf 6} (1985).

\bibitem{krieger} {\sc W. Krieger.} On quasi-invariant measures in uniquely 
ergodic systems. {\em Invent. Math.} {\bf 14} (1971), 184-196.

\bibitem{KS} {\sc J. Kwiatkowski \& A. Siemaszko.} Discrete orbits for topologically 
transitive cylindrical transformations. {\em Discrete and Contin. Dyn. Syst.} {\bf 27} 
(2010), 945-961. 

\bibitem{le-calvez-yoccoz} {\sc P. Le Calvez \& J.-C. Yoccoz.} Un th\'eor\`eme d'indice pour 
les hom\'eomorphismes du plan au voisinage d'un point fixe. {\em Ann. of Math.} {\bf 146} 
(1997), 241-293.

\bibitem{MS} {\sc S. Matsumoto \& M. Shishikura.} Minimal sets of certain annular 
homeomorphisms. {\em Hiroshima Math. J.} {\bf 32} (2002), 207-215.

\bibitem{OP} {\sc J. Moulin-Ollagnier \& D. Pinchon.} Syst\`emes dynamiques topologiques 
I. \'Etude des limites de cobords. {\em Bulletin de la S.M.F.} {\bf 105} (1977), 405-414.

\bibitem{book} {\sc A. Navas.} {\em Groups of Circle Diffeomorphisms.} 
Chicago Lectures in Mathematics (2011).

\bibitem{OO} {\sc G. Ochs \& V.I. Oseledets.} Topological fixed point theorems do not 
hold for random dynamical systems. {\em J. Dynam. Differential Equations} {\bf 11} 
(1999), 583-593.

\bibitem{perez-marco} {\sc R. P\'erez-Marco.} Fixed points and circle maps. 
{\em Acta Math.} {\bf 179}, number {\bf 2} (1997), 243-294.

\bibitem{ponce} {\sc M. Ponce.} Local dynamics for fibred holomorphic transformations. 
{\em Nonlinearity} {\bf 20}, number {\bf 12} (2007), 2939-2955.

\bibitem{sch} {\sc L.G. Shnirelman.} An example of a transformation of the plane 
(Russian), Proc, Don. Polytechnic Inst. (Novochekassk) {\bf 14}, Science section, 
Fis-math. part (1930), 64-74.

\bibitem{yoccoz} {\sc J.-C. Yoccoz.} Sur la disparition de propri\'et\'es de type 
Denjoy-Koksma en dimension $2$. {\em C. R. Acad. Sci. Paris} {\bf 13} (1980), 655-658.

\end{thebibliography}
\end{document}